\documentstyle[12pt,a4]{article}

\def\qed{\hspace*{0pt}\hfill
\hbox{\begin{picture}(10,10)\put(0,0){\framebox(6,6){}}
\end{picture}}\medskip}

\begin{document}

\setlength{\parskip}{0.5cm}

\begin{center} {\bf
BOLZANO-WEIERSTRASS PRINCIPLE OF CHOICE EXTENDED TOWARDS ORDINALS}
\vspace{0.5cm}

 by

\vspace{0.5cm}

{\it        W\l adys\l aw Kulpa, Szymon Plewik} and {\it Marian
Turza\'nski } \end{center}

\vspace{0.5cm}

{\bf Abstract.} The Bolzano-Weierstrass principle of choice is the oldest method of the set theory, traditionally used in mathematical analysis. We are extending it towards transfinite sequences of steps indexed by ordinals. We are introducing the notions: hiker's tracks, hiker's maps and statements $P_n(X, Y, m)$; which are  used similarly in finite, countable and uncountable cases. New proofs of Ramsey's theorem and Erd\"{o}s-Rado theorem are presented as some applications.

\vspace{0.5cm}

{\bf  I - Introduction}. The Bolzano-Weierstrass principle of choice is based on succeeding divisions of a segment onto disjoint subsegments and on choice of a subsegment which has some desired quality. Our extension towards transfinite sequences of steps indexed by ordinals imitates a hiker's track. Any hiker's step corresponds to dividing. It is uniquely determined by his
previous steps and by his destination point. To express our extension we
consider statements $P_{n}(X, Y, m )$. If some of these statements depend on uncountable parameters,
then we use them for a new proof of the Erd\"os-Rado partition
theorem.  If some others depend on countable parameters, then we 
 use them for a proof of Ramsey's theorem. But, if any one depends on finite
parameters, then we could introduce numbers $p(k,r,n)$ which are very similar to the so called Ramsey's numbers.

 For a given set $X$  denote its cardinality by $|X|$. If $n$ is  a
natural number, then   $[X]^1=X$, and $ [X]^n = \{ b \subseteq X:
|b| = n \} $. Infinite ordinals are usually denoted by Greek letters. Sometimes we write $\alpha \in \beta$, instead of $\alpha < \beta$. The remaining notations are standard.

\vspace{0.5cm}

{\bf  II - Hiker's track, hiker's map.}
Let $X$ be a well ordered set, i.e. $X$ is an ordinal number. 
Fix a function $\displaystyle f:[X]^{n+1} \rightarrow Y$. 
Let $ x_0 <  x_1 < \ldots <  x_{n-1} $ be the first $n$ points in $X$. For any point 
$x \in X$ the increasing sequence $$\displaystyle \{ x_{\beta}: \beta \leq \delta (x) \} \subseteq X $$ 
is called the $ x$-\textit{hiker's track}, if any $\displaystyle x_{\beta}$ where $\beta \geq n$,  is defined us  follows.
Suppose that $\displaystyle \{ x_{\gamma}: \gamma < \beta \}$ has been defined. 
For any $\displaystyle s \in [\beta]^n$ consider the subset 
$$ \delta (x, s) =  \left\{ y \in X: \quad f( \{
x_{\gamma}: \gamma \in s \} \cup \{ y \} ) = f  ( \{ x_{\gamma}: 
\gamma \in s \} \cup  \{ x \}) \right\} \subset X. $$ 
Let $x_\beta$ be the first point in the intersection
$ \cap \{ \delta (x, s): s \in [\beta]^n \} \subset X.$ Our construction stops when $x_{\beta} = x$. This $\beta$ is denoted $\delta (x)$.
Any $ x$-hiker's track is uniquely
determined by the 
increasing sequence  $ \{ x_{\gamma}:
\gamma \leq   \delta (x )  \} \subseteq X $ and $
x_{\delta (x)} = x$ holds.

 If $s \in [\delta (x)]^n$, then
put $\displaystyle f_{x} (s) = f( \{ x_{\zeta}: \zeta \in s\} \cup \{ x \}
).$ A function $f_x$ is called $x$-\textit{hiker's map}.

{\bf Theorem 1}. {\it Any hiker's map $\displaystyle f_{x}:  [\delta (x ) ]^n \to Y  $ is uniquely determined by the subsequence $ \{ x_{\gamma}:  \gamma < \delta (x)  \}$ of the $ x$-hiker's track. In other words, the map which for each point $x\in X$ assigns $f_x$ is one-to-one.}

{\bf Proof}.  We shall prove that 
 if $x \not= y$, then functions $f_x$ and $f_y$ are different. Indeed, if $\delta (x) \not= \delta(y)$, then the functions have different domains. If $\delta (x) = \delta(y)$, then the $ x$-hiker's track $\displaystyle \{ x_{\beta}: \beta \leq \delta (x) \}$ and the $ y$-hiker's track
$\displaystyle \{ y_{\beta}: \beta \leq \delta (y)\}$ are different. Suppose $\beta$ is the first ordinal such that $x_\beta \not= y_\beta$. Without the loss of generality, assume $x_\beta <y_\beta$. Thus there exists $s\in [\beta]^n$ such that $\displaystyle x_\beta \not\in \delta (y, s)$. Hence  
$$f( \{ y_{\zeta}: \zeta \in s\} \cup \{ x_\beta \}
) \not= f( \{ y_{\zeta}: \zeta \in s\} \cup \{ y \}
). $$
But 
$$ f( \{ x_{\zeta}: \zeta \in t\} \cup \{ x \}
)=f( \{ x_{\zeta}: \zeta \in t\} \cup \{ x_\beta \}
),$$ for each $t\in [\beta]^n$. 
In consequence,
$$f( \{ y_{\zeta}: \zeta \in s\} \cup \{ y \}
) \not= f( \{ x_{\zeta}: \zeta \in s\} \cup \{ x \}
). $$
For such $s$ the following 
$\displaystyle f_x (s ) \not= f_y (s ) $ holds.
In others words, the
  hiker's maps $f_x$ and $f_y$ are different. \qed 

 In the literature, when $X$ is a countable set, the above construction is called the  Bolzano-Weierstrass principle of choice. In this way  the method of indicating a monotone sequence which is contained in an infinite set of real numbers is honored. In this note we are extending this principle onto an arbitrary ordinal number instead of a countable set $X$.

\hspace{0.5cm}

{\bf III - Statements $P_{n}(X, Y, m )$.}
Fix some sets $X$, $Y$ and $m$. Assume that the set $X$ is well ordered and $m$ is an ordinal number. For any natural number $n$ consider 
the following statement.

 {\it For each function $f:[X ]^{n+1} \rightarrow Y $ there exists an increasing sequence $\{
x_{\beta}: \beta < m \} \subseteq X$ such that if $\beta_0 < \beta_1
< \ldots < \beta_n < \beta_{n+1} < m $, then $$ f( \{
x_{\beta_{0}}, x_{\beta_{1}}, \ldots , x_{\beta_{n-1}}, x_{\beta_{n}}
\} ) =f(\{ x_{\beta_{0}}, x_{\beta_{1}}, \ldots , x_{\beta_{n-1}},
x_{\beta_{n+1}} \}) .$$ } Denote this statement by
 \textbf{$P_{n}(X, Y, m)$}.
In other words, 
 \textbf{$P_{n}(X, Y, m)$} means that for each function $f:[X ]^{n+1} \rightarrow Y $ there exists a hiker's track which as a sequence has the length $m +1$. Note that in this definition $X$ does not need to be well ordered. But in applications we usually assume that $X$ is well ordered. 
 
 If $n=0$, then one has the  pigeonhole principle. Then $ P_0(m^+, m,m^+ )$ holds for any infinite cardinal number $m$. 
Because of $m<m^+$, then 
 for any function 
$f:m^+\rightarrow m $ there exists a point $x \in m$ such that the preimage $f^{-1} (x)$ has the cardinality $m^+$. For a such $x$ one can write $f^{-1} (x)= \{ a_{\alpha}: \alpha < m^+\} $, i. e.
$x= f(a_\beta) = f(a_\alpha)$, whenever $\alpha < \beta < m^+$.

{\bf  Theorem 2}.
 {\it If    $n>0$ and $\lambda$ is an infinite cardinal
number, then the statement $P_{n}( |2^{\lambda}|^+, 2^{\lambda},
\lambda^+ + 1)$ holds}.

{\bf Proof.} Fix a function $f:[|2^{\lambda}|^+ ]^{n+1} \rightarrow 2^{\lambda}$ and put  $X = |2^{\lambda}|^+ $ and $Y= 2^{\lambda}$ and
$ x_0 =0,\quad  x_1 =1,\quad  \ldots , \quad x_{n-1} = n-1.$ 
For any ordinal number  
$\alpha \in |2^{\lambda}|^+$ consider $\alpha$-hiker's track
$\displaystyle \{ a_{\beta}: \beta \leq \delta (\alpha) \}$. 
 Suppose that $|\delta(\alpha)| \leq \lambda$ for any $\alpha \in |2^{\lambda}|^+$. There are at most $\lambda^+$ different ordinals of the form $\delta(\alpha)$ and at most $2^\lambda$ different hiker's maps $f_{\alpha}:  [\delta (\alpha ) ]^n \to 2^\lambda  $. Because of $\lambda^+ \cdot 2^\lambda = 2^\lambda$ one obtains a contradiction, since Theorem 1 says that there has to be $(2^\lambda)^+$ different hiker's maps. So, there exists $\alpha \in |2^\lambda|^+$ such 
that $|\delta(\alpha)| = \lambda^+$.  \qed

{\bf  Theorem 3}.
 {\it Let    $n>0$ and $\lambda$ be an infinite cardinal
number. If $m$ is a cardinal number such that $2^m \leq 2^{\lambda}$, then   $P_{n}( |2^{\lambda}|^+, 2^{\lambda},
m+1)$ holds}.

{\bf Proof.} Fix a function $f:[|2^{\lambda}|^+ ]^{n+1} \rightarrow 2^{\lambda}$ and put  $X = |2^{\lambda}|^+ $ and $Y= 2^{\lambda}$ and
$ x_0 =0,\quad  x_1 =1,\quad  \ldots , \quad x_{n-1} = n-1.$ 
For any ordinal number  
$\alpha \in |2^{\lambda}|^+$ consider $\alpha$-hiker's track
$\displaystyle \{ a_{\beta}: \beta \leq \delta (\alpha) \}$. 
 Suppose that $|\delta(\alpha)| < m$ for any $\alpha \in |2^{\lambda}|^+$. There are at most $m$ different ordinals of the form $\delta(\alpha)$ and at most $2^\lambda = (2^\lambda)^m$ different functions $f_{\alpha}:  [\delta (\alpha ) ]^n \to 2^\lambda  $. Because of $m \cdot 2^\lambda = 2^\lambda$ one obtains a contradiction, since by Theorem 1, there are $(2^\lambda)^+$ different hiker's maps. So, there exists $\alpha < |2^\lambda|^+$ such 
that $|\delta(\alpha)| = m$. \qed

{\bf  Theorem 4}.
 {\it Let    $n>0$ and $\lambda$ be an infinite cardinal
number. If $2^m \leq 2^\lambda$ for any cardinal number $m < 2^{\lambda}$, then $P_{n}( |2^{\lambda}|^+, 2^{\lambda},
2^\lambda +1)$ holds}.

{\bf Proof.} Fix a function $f:[|2^{\lambda}|^+ ]^{n+1} \rightarrow 2^{\lambda}$ and put  $X = |2^{\lambda}|^+ $ and $Y= 2^{\lambda}$ and
$ x_0 =0,\quad  x_1 =1,\quad  \ldots , \quad x_{n-1} = n-1.$ 
For any ordinal number  
$\alpha \in |2^{\lambda}|^+$ consider $\alpha$-hiker's track
$\displaystyle \{ a_{\beta}: \beta \leq \delta (\alpha) \}$. 
 Suppose that $|\delta(\alpha)| < 2^\lambda$ for any $\alpha \in |2^{\lambda}|^+$. 
There are at most $2^\lambda$ different ordinals of the form $\delta(\alpha)$ and at most $$(2^\lambda)^{|\delta(\alpha)|} = 2^{\lambda \times |\delta(\alpha)|} = 2^{|\delta(\alpha)|\times \lambda } = (2^\lambda)^\lambda= 2^\lambda$$ different functions $f_{\alpha}:  [\delta (\alpha ) ]^n \to 2^\lambda  $. One obtains another contradiction with Theorem 1. This  follows that there exists $\alpha < |2^\lambda|^+$ such 
that $|\delta(\alpha)| = 2^\lambda$. \qed

\vspace{0.5cm}

{\bf  IV - Applications to some  proofs of  Erd\"os--Rado partition theorems}.
 To give some applications of statements $P_{n}(X, Y, m )$ we start with a proof of P. Erd\"os and R. Rado theorem [2]. Let $\exp^{(0)}( X) = X $ and $\exp^{(n+1)}(X ) =
2^{\exp^{(n)}( X)}$. We need statements $$P_{k}( |\exp^{(k)}(\zeta )|^+ , \zeta, |\exp^{(k-1)}(\zeta)|^+
),$$ where  $0<k \leq n$. One can deduce these statements from  Theorem 2
putting $\lambda = \exp^{(k-1)} (\zeta ) $ and restricting the second 
parameter to $\zeta < | 2^{\lambda}|$. 

{\bf  Theorem} (P. Erd\"os and R. Rado [2]).
 {\it Let $n$ be a natural number, but $ \zeta $ and $\kappa$ be infinite
cardinal numbers, and assume that $\kappa  > |\exp^{(n)} (\zeta)| .$ Then for
any  function $ f:[\kappa ]^{n+1} \to \zeta  $ there exist an ordinal
$\varphi < \zeta $ and a subset $ Z \subseteq \kappa $ such that   $|Z| >
\zeta $ and $ [Z]^{n+1} \subseteq f^{-1}(\varphi).$}

{\bf Proof.} We proceed by induction on $n$. For $n=0$ we have assumed
$\kappa > \zeta$ and $P_0(\zeta^+, \zeta, \zeta^+)$ holds by the same argumentat as this before Theorem 2. Let
$|\exp^{(n-1)}(\zeta)| =\lambda$, and assume that the theorem holds for a
natural number $n-1\geq 0$.  Since  $ |2^{\lambda}|^+ \leq \kappa $ the statement
 $P_{n}( |2^{\lambda}|^+ , \zeta, \lambda^+ )$ yields a hiker's track
$\{ a_{\beta}: \beta < \lambda^+ \} $ such that if $\beta_0 < \beta_1 <
\ldots < \beta_n < \beta_{n+1} < \lambda^+ ,$ then  $$ f( \{
a_{\beta_{0}},a_{\beta_{1}}, \ldots , a_{\beta_{n-1}}, \alpha_{\beta_{n}}
\} ) =f(\{ a_{\beta_{0}}, a_{\beta_{1}}, \ldots , a_{\beta_{n-1}},
a_{\beta_{n+1}} \}) .$$ Consider the notion of a hiker's map, i.e. if $\beta_0 < \beta_1 < \ldots  <
\beta_{r-1} < \beta < \lambda^+ $, then  put  $$ F( \{ \beta_{0},\beta_{1}, \ldots ,
\beta_{n-1} \} ) =f(\{ a_{\beta_{0}}, a_{\beta_{1}}, \ldots ,
a_{\beta_{n-1}},  a_{\beta} \}) . $$ This hiker's map is a function $F: [
\lambda ^+ ]^n \to \zeta $. By the induction hypothesis there exist
 an ordinal $\varphi < \zeta $ and a subset $S \subseteq  \lambda^+ $ such
that $[S]^n \subseteq F^{-1} (\varphi) $   and $|S|>\zeta$.  The subset
$ Z = \{ \alpha _{\beta}:  \beta \in S \} \subseteq \kappa$ has
cardinality greater than $\zeta$, and if $ \{\beta_0, \beta_1, \ldots,
{\beta}_{k-1}, \beta \} \subseteq S $ and $ {\beta}_0 < \beta_1 < \ldots <
{\beta}_{k-1} < \beta $, then $$f ( \{ \alpha_{\beta_{0}},
\alpha_{\beta_{1}}, \ldots , \alpha_{\beta_{n-1}}, \alpha_{\beta} \}) = F
( \{ \beta_{0}, \beta_{1}, \ldots , \beta_{n-1}
 \}  ) = \varphi.$$ This clearly implies $[Z]^{n+1} \subseteq
f^{-1}(\varphi)$. \qed

In the literature there are many proofs of the  of the Erd\"os--Rado
partition theorems. Our proof  looks most  similar to that  of 
J.D.Monk, see [3] p. 1230. This similarity could be understood such 
that  the notion of "pre-homogeneous" is replaced by suitable statements 
$P_n( X , Y, \lambda )$. This  gives some reasons to consider $P_n(
X , Y, \lambda )$ as self-made notions.

Add two corollaries which are similar to some results which are in the book [1], pages 8 - 11.

{\bf  Corollary 5}. 
	  {\it Let  $m$, $\xi$, $ \lambda $ and $\kappa$ be infinite
cardinal numbers. Assume that $\kappa  > |2^\lambda|^+ $ and $\xi<m$ and $2^m \leq 2^\lambda$. Then for
any  function $ f:[\kappa ]^2 \to \xi  $ there exist an ordinal
$\varphi < \xi $ and a subset $ Z \subseteq \kappa $ such that $|Z| \geq
m $ and $ [Z]^{2} \subseteq f^{-1}(\varphi).$}

 \textbf{Proof}. By Theorem 3 the statemnet $P_1( \kappa, \xi, m+1)$ holds. This yields a hiker's track 
 $\{ a_{\beta}: \beta \leq m \} $ such that if $\beta < \gamma <
m ,$ then  $$ f( \{
a_{\beta},a_{\gamma}
\} ) =f(\{ a_{\beta}, a_{m}\}) .$$ In consequence one obtains a hiker's map $F: m \to \xi $, where $ 
F(  \beta  ) =f(\{ a_{\beta}, a_{m})$ for any $\beta <m $. By the pigeonhole principle and since $\xi < m$ 
there exist
 an ordinal $\varphi < \xi $ and a subset $Z \subseteq  \{ a_{\beta}: \beta \leq m \} $ such
that $[Z]^2 \subseteq f^{-1} (\varphi) $ and $|Z|= m$. \qed

{\bf  Corollary 6}.  
	{\it Let   $\xi$, $ \lambda $ and $\kappa$ be infinite
cardinal numbers. Assume that $\kappa  > |2^\lambda|^+ $ and $\xi<2^\lambda$ and  $m < 2^\lambda$ always implies that $2^m \leq 2^\lambda$. Then for
any  function $ f:[\kappa ]^2 \to \xi  $ there exist an ordinal
$\varphi < \xi $ and a subset $ Z \subseteq \kappa $ such that   $|Z| \geq |2^\lambda| $ and $ [Z]^{2} \subseteq f^{-1}(\varphi).$}

  \textbf{Proof}. By Theorem 4  the statement  $P_1( \kappa, \xi, 2^\lambda +1)$ holds. This yields a hiker's track 
 $\{ a_{\beta}: \beta \leq 2^\lambda \} $ such that if $\beta < \gamma <
2^\lambda ,$ then $$ f( \{
a_{\beta},a_{\gamma}
\} ) =f(\{ a_{\beta}, a_{2^\lambda}\}) .$$ In consequence one obtains a hiker's map $F: 2^\lambda \to \xi $, where $ 
F(  \beta  ) =f(\{ a_{\beta}, a_{2^\lambda}\})$ for any $\beta <2^\lambda $. By the pigeonhole principle and since $\xi < 2^\lambda$   
there exist
 an ordinal $\varphi < \xi $ and a subset $Z \subseteq  \{ a_{\beta}: \beta \leq 2^\lambda \} $ such
that $[Z]^2 \subseteq f^{-1} (\varphi) $ and $|Z|= 2^\lambda$. \qed

{\bf V - On a proof of Ramsey's theorem.} In this part we give applications of statements $P_n(\omega , r , \omega )$, where $\omega$ denotes the set natural numbers and $r$ is a natural number. To do this we present a proof of the Ramsey's theorem, see [5]. In [4] p. 5 there is a proof of
Ramsey's theorem which contains some aspects of the Bolzano--Weierstrass
principle of choice. 

{\bf Theorem 7}.
{\it If $r>0$ and  $n$ are natural numbers but  $\omega$ is the first
infinite ordinal, then the statement $P_n(\omega , r , \omega)$ holds.}

{\bf Proof.} Fix a function $f:[\omega]^{n+1} \to r$. Any  natural number
$k \geq n$
 uniquely determines the hiker's map  $f_{k}:  [\delta (k ) ]^n
\to r $. All maps  $f_{k} $ forms an infinite tree. By the
K\"onig infinity lemma  this tree possesses  an infinite
 path. Any such infinite path marks a desired hiker's track. \qed

 Now, using the reduction from our proof of the
Erd\"os--Rado partition theorem: the reduction from $P_n(\omega , r , \omega )$
to $P_{n-1}(\omega , r , \omega )$); we obtain a proof of the Ramsey's theorem.

\textbf{Ramsey's theorem}. (F. P. Ramsey [5]). {\it If $r>0$ and $n$ are natural numbers, then for any function $ f:[\omega
]^{n+1} \to r  $ there exist a natural number $m $ and an infinite subset
$ Z \subseteq \omega $ such that $[Z]^{r+1} \subseteq f^{-1}(m).$ }

\textbf{Proof}. 
We proceed  by induction on $n$. For $n=0$ we have the pigeonhole principle. 
Assume that Ramsey's  theorem holds for a
natural number $n-1\geq 0$.  The statement
 $P_{n}( \omega , r, \omega)$ yields an infinite hiker's track
$\{ a_{k}: k < \omega \} $ such that if $k_0 < k_1 <
\ldots < k_n < k_{n+1} ,$ then  $$ f( \{
a_{k_{0}},a_{k_{1}}, \ldots , a_{k_{n-1}}, \alpha_{k_{n}}
\} ) =f(\{ a_{k_{0}}, a_{k_{1}}, \ldots , a_{k_{n-1}},
a_{k_{n+1}} \}) .$$ To the  hiker's map $F: [
\omega ]^r \to n $ , where $$ F( \{ k_{0},k_{1}, \ldots ,
k_{n-1} \} ) =f(\{ a_{k_{0}}, a_{k_{1}}, \ldots ,
a_{k_{n-1}},  a_{k_{n}} \})  $$ one applies  the induction hypothesis. \qed

{\bf Corollary 8}.
{\it If $r>0$ and  $n$ are natural numbers but  $\omega$ is the first
infinite ordinal, then the statement $P_n(\omega , r , \alpha )$
holds for any countable ordinal \mbox{number $\alpha$. \hfill \qed}}

{\bf VI -  Numbers $p(k,r,n)$.}   Consider  statements $P_n(X,  Y,
\lambda )$ for cases when $X$, $Y$ and $\lambda$ are natural numbers. Fix
positive natural numbers $k$, $r$ and $n$. Similar to the definition of
Ramsey's numbers - compare [4] p.  13 - let $p(k,r,n)>r $ be the least natural number such that the statement  $P_r(p(k,r,n),n,k)$  holds.  This means that $p(k,r,n)$ is the least natural number such that for
any function $f:[p(k,r,n) ]^{r+1} \to n $ there exists a hiker's track $\{
a_{i}: i < k \}  $ such that  $$ f( \{
a_{i_{0}},a_{i_{1}}, \ldots , a_{i_{r-1}}, \alpha_{i_{r}}
\} ) =f(\{ a_{i_{0}}, a_{i_{1}}, \ldots , a_{i_{r-1}},
a_{i_{r+1}} \}) ,$$ whenever  $i_0 < i_1 < \ldots < i_r
< i_{r+1}  .$  Numbers $p(k,r,n)$ are well defined since 
the following holds.

{\bf Theorem 9}.
{\it If $n$, $r$ and $k$ are positive natural numbers, then
$$p(k,r,n) - r < n^{r+0\choose {r}} +n^{r+1\choose {r}} + \ldots
+n^{r+k-2\choose {r}}+1 .$$}

{\bf Proof.} Use again the Bolzano-Weierstrass principle of choice.  If $r\leq k $,
then the function   $f_{k}:  [\delta (k ) ]^r \to n  $  is
uniquely determined. Also, if $\delta (k ) =r+i$, then  there are
$n^{r+i\choose {r}} $ possibilities for   any $f_{k}$ with the
domain of cardinality $i$. Therefore $  p(k,r,n) -r -1 \geq
n^{r+0\choose {r}} + n^{r+1\choose {r}} + \ldots +n^{r+k-2\choose {r}}  $
implies that a function $f_{k}$ with the
domain of cardinality $k-1$ has to be defined . Any such a function $f_{k}$ designs a
desired sequence. \qed

{\bf References.}

[1] W. W. Comfort and S. Negrepointis, {\it Chain conditions in topology }, Cambridge University Press (1982).

[2] P.  Erd\"os and R. Rado, {\it A partition calculus in set theory},
Bull. Amer. Math. Soc. 62(1956), 427-489.

[3]  J. D. Monk, {\it Appendix on set theory}, in: Handbook of Boolean
algebras, Elsevier Science Publishers  (1989), 1213 -1233.

[4] H.  J. Pr\"omel and B. Voigt,  {\it Aspects of Ramsey--theory I: sets
}, Forschungsinstitut f\"ur Diskrete Mathematik Institut f\"ur
\"Okonometrie und Operations Research Rheinische
Friedrich--Wilhelms--Universit\"at Bonn, report No 87 495--OR (1989).

[5] F. P. Ramsey, \textit{On a problem of formal logic}, Proc. London Math. Soc. 2 (1930), 264 - 286.

\begin{tabular}{lcl}
Department of Mathematics & \mbox{} \hfill  \mbox{}&  e-mail
   addresses: \\ Silesian University && kulpa@ux2.math.us.edu.pl \\  ul.
 Bankowa 14 &&        plewik@ux2.math.us.edu.pl \\ 40-007 Katowice,
Poland & & mtturz@ux2.math.us.edu.pl \\ &&\\ &&\\ AMS Subject
Classification (1991):  &&Primary:   05E20; \\ && Secondary: 03E05,
04A20.  \end{tabular}

\end{document}